\documentclass{article}
\usepackage[usenames,dvipsnames]{color}
\usepackage{graphicx}
\usepackage{camnum}
\usepackage{harvard}
\usepackage{amsmath}
\usepackage{eufrak}
\usepackage{soul}
\usepackage{booktabs}
    \def\O#1{{\cal O}\left(#1\right)}

\newcommand{\ee}{{\mathrm e}}
\newcommand{\ii}{{\mathrm i}}
\newcommand{\ad}{\mathrm{ad}}

\newcommand{\dexp}{\mathrm{dexp}}

\newcommand{\sBCH}{\CC{sBCH}}
\newcommand{\A}{\mathcal{A}}
\newcommand{\K}{\mathcal{K}}

\newcommand{\DDD}{\mathcal{D}}
\newcommand{\WW}{\mathcal{W}}
\newcommand{\dx}{\partial_x}
\newcommand{\ve}{\varepsilon}
\newcommand{\schr}{Schr\"odinger }

\newcommand{\dummy}[1]{#1}
\newcommand{\PSc}[1]{\textcolor{black}{#1}}
\newcommand{\KKc}[1]{\textcolor{black}{#1}}
\newcommand{\PBc}[1]{\textcolor{black}{#1}}

\newcommand{\Int}[4]{\int_{#2}^{#3}\!#4\,\mathrm{d}#1}
\newcommand{\Ang}[2]{\left\langle #2 \right\rangle_{#1}}

\newcommand{\ang}[2]{\left\langle #1 \right\rangle_{#2}}

\begin{document}
\title{Efficient methods for time-dependence in semiclassical Schr\"odinger equations}

\author{Philipp Bader,\footnote{La Trobe University, Department of Mathematics, Kingsbury Dr, Melbourne 3086 VIC, Australia}\ \  Arieh Iserles,\footnote{Department of Applied Mathematics and Theoretical Physics, University of Cambridge, Wilberforce Rd, Cambridge CB3 0WA, UK.}\ \  Karolina Kropielnicka\footnote{Institute of Mathematics, University of Gda\'nsk, Wit Stwosz Str.\ 57, 90-952 Gda\'nsk, Poland.}\ \  \& Pranav Singh\footnote{Department of Applied Mathematics and Theoretical Physics, University of Cambridge, Wilberforce Rd, Cambridge CB3 0WA, UK.}}

\maketitle

\begin{abstract}
	We build an efficient and unitary (hence stable) method for the solution of the semi-classical Schr\"odinger equation subject with explicitly time-dependent potentials.
	The method is based on a combination of the Zassenhaus decomposition \cite{bader14eaf} with the Magnus expansion of the time-dependent Hamiltonian.
	We conclude with numerical experiments.
\end{abstract}

\section{Introduction}
Rapid advances in laser technologies over the recent years have led to a significant progress in the control of systems at the molecular level \cite{shapiro03qc}.
\PSc{Pioneering work in the control of chemical systems at the quantum level was done in the study of photo-dissociation and bimolecular reactions.
Various control techniques such as the pump-dump quantum control scheme \cite{Kosloff1989201} and the coherent control schemes \cite{Shapiro} have had numerous experimental validations and applications \cite{Zhu77,Vogt2006211}.}

\PSc{These experimental successes and a dramatic improvement in our ability to shape femtosecond laser pulses over the recent years has led to a great deal of interest in the development of a systematic way of designing controls (shaping laser pulses) and a requirement for rigourous mathematical analysis of issues such as controllability \cite{lebris}.}

%of numerical schemes for a systematic procedure of design of controls through techniques such as optimal control are becoming popular
%
% and the potential applicability to a vast range of quantum systems in quantum chemistry, solid state physics and quantum computing
%
%numerical schemes for systematically and the various mathematical challenges that the problem of optimal control throws up

\PSc{In the case of laser induced breakdown (photo-dissociation) of a molecule, for instance, there is a great deal of interest in designing lasers that achieve efficient breakdown. The fact that the dissociation timescales are often themselves in femtoseconds means that it cannot generally be assumed that the laser pulse causes near-instantaneous and efficient excitation of a molecule sitting in the ground state, having no other influence thereafter -- the correct dynamics require taking into account the time-dependent nature of the electric potential (laser) throughout the evolution of the wavefunction.}

\PSc{To analyse the control exerted by these lasers we need efficient means of computing the \schr equation featuring time-dependent Hamiltonians, existing strategies for which are either low accuracy or become prohibitively expensive with higher orders of accuracy. }

\PSc{Optimal control schemes for designing laser pulses is often posed as an inverse problem that is solved via optimisation schemes requiring repeated solutions of \schr equations with modified time-dependent Hamitlonians. An ability to efficiently solve these \schr equations with moderately large time steps and high accuracy becomes crucial here, creating a need for high-order methods \cite{karlsson}. }

\PSc{In this paper,} we are interested in the numerical computation of the linear, time-dependent Schr\"o\-din\-ger equation in a semiclassical regime for a nucleus moving in a time-dependent electric field,
\begin{equation}\label{eq.se}
\partial_t u(x,t) = -\frac{\ii}{\varepsilon} H u(x,t)= \ii[\varepsilon \Delta - \varepsilon^{-1} V(x,t) ]u(x,t),\ x\in[-1,1],\ t\geq0,
\end{equation}
equipped with an initial condition $u(x,\KKc{0})=u_0(x)$, and periodic boundary conditions. We assume, that the potential $V\KKc{(\cdot,t)}\in C^{\infty}[-1,1]$ is periodic\KKc{, for $t\geq0$}.

The equation (\ref{eq.se}) is posed on a Hilbert space $\mathcal{H} = \CC{L}_2[-1,1]$, and the squared modulus of the solution is the probability density of finding the particle in state $x$ at time $t$. For this reason, the initial \KKc{condition $u_0(x)$} %density
is normalised to one and it is easy to see that the norm of the solution is an invariant,
$$
\|u(x,t)\|^2_{\CC{L}^2}=\int_{-1}^{1}|u(x,t)|^2dx=\|u(x,0)\|^2_{\CC{L}^2}.
$$
The wave function undergoes unitary evolution, which we wish to preserve under discretisation -- both because of physical significance, and since, as we mention in Section~4, it implies stability of the numerical method.

\PSc{Here the semiclassical parameter $0<\varepsilon\ll1$ may arise out of a {\em Born--Oppenheimer approximation}, via spatio-temporal scaling, or a combination of the two, depending upon the physical system under consideration. }
The regularity of $V$ depends on the order of desired accuracy, but for convenience we have assumed that it is smooth in its domain.
The initial condition is usually a high-frequency wave packet, but even if it is non-oscillatory it can be shown, cf. the analysis in \cite{jin11mac}, that the solution to this Schr\"odinger equation is highly oscillatory, with frequency of at least $\mathcal{O}(\varepsilon^{-1})$. This, as a matter of fact, is the main reason why finding an effective numerical method for (\ref{eq.se}) is such a challenging task. Obviously, the naive approach of finite differences is out of consideration, but instead the usual methodology consists of a semidiscretisation in space \PSc{via spectral methods} followed by an exponential splitting.

\KKc{The first step in approximating (\ref{eq.se}) usually is spatial discretisation, which yields the following system of ODEs
\begin{equation}\label{eq:SD}
 \MM{u}'(t)=\ii (\ve\K^2 - \ve^{-1} \DDD_{V(\cdot,t)})\MM{u}(t),\qquad t\geq0,
\end{equation}
where  $\K^2$ and $\DDD_{V(\cdot,t)}$ are $M \times M$ matrices representing the discretisation of second derivative and the  multiplication by $V(\cdot,t)$, respectively. We understand, that $\MM{u}(t)\in\BB{C}^M$ is a vector representing an approximation to the solution (\ref{eq.se}) at time $t$ and $\MM{u}(0)$ is derived from the initial conditions.
%At that stage we have to possibilities: either to opt for higher order method, which is computationally costly, or for lower order method, but much cheaper.
A second order method can be obtained by freezing the matrix $\DDD_{V(\cdot,t)}$ in the middle of interval $[0,t]$ and applying the Strang splitting
%$$
%\ee^{t \ve\K^2-t\ve^{-1} \DDD_{V(\cdot,t)}}=\ee^{t\ve\K^2}\ee^{-t\ve^{-1} \DDD_{V(\cdot,t/2)}}+\mathcal{O}(t^2),
%$$
%or second order Strang splitting
$$
%\ee^{t\ve\K^2-t\ve^{-1} \DDD_{V(\cdot,t)}}
\mathbf{u}(t)=\ee^{\frac12t\ve\K^2}\ee^{-t\ve^{-1} \DDD_{V(\cdot,t/2)}}\ee^{\frac12 t\ve\K^2}\PSc{\MM{u}(0)}+\mathcal{O}(t^3).
$$
This splitting has the advantage of separating scales ($\ve$ and $\ve^{-1}$) as well as easily computable exponentials. Using spectral collocation or spectral spatial discretisation methods, the matrices $\K$ and $\DDD_{V(\cdot,t/2)}$ are either diagonal (thus exponentiated directly), or circulant (thus approximated by FFT).
Higher order methods can be obtained by representing the solution to (\ref{eq:SD}) through a Magnus expansion,
\begin{displaymath}
 \MM{u}(t)=\ee^{\MM{\Theta}(t)}\MM{u}(0),
\end{displaymath}
where $\MM{\Theta}(t) \in \mathfrak{u}_M(\BB{C})$ is a time-dependent $M \times M$ skew-Hermitian matrix obtained as an infinite series $\sum_{k=1}^\infty \MM{\Theta}^{[k]}(t)$ with each $\MM{\Theta}^{[k]}(t)$ composed of $k$ nested integrals and commutators of the matrices $\ii\ve\K^2$ and $\ii\varepsilon^{-1} \DDD_{V}$. }

\KKc{This approach was exploited in \citeasnoun{hochbruck03omi}, where authors conclude that the Magnus expansion $\MM{\Theta}(t) $ is convergent if $t$ is such, that for some constant $c$ the inequality $t \| \K \| \leq c$ holds.}
\KKc{This, as explained later, forces us to time step of order $\mathcal{O}(\ve)$.}
\KKc{Another serious drawback of this approach lies in the costly approximation of the exponential $\ee^{\MM{\Theta}(t)}$. As it occurs,  the  exponent $\MM{\Theta}(t)$ ends up to be of a large size (both: spectral and dimensional), and neither diagonal nor circulant. Indeed, observe, that the highly oscillatory nature of the solution to (\ref{eq.se}) requires a large number of degrees of freedom in the spatial discretisation, $M = \O{\ve^{-1}}$. Since the differentiation matrix $\K$ scales as $\O{M}=\O{\ve^{-1}}$, the operator $\MM{\Theta}(t)$, as a sum of nested commutators of $\ii\ve\K^2$ and $\ii\varepsilon^{-1} \DDD_{V(\cdot,t)}$, occurs to be  a large matrix which does not posses any favourable structure that could allow an effective approximation of the exponential $\exp(\MM{\Theta}(t))$.}

Powerful tools like Zassenhaus splitting or Baker--Campbell--Hausdorff formula were historically avoided in splitting methods due to the large computational cost of nested commutators. However as it happens, choosing the correct, infinite--dimensional Lie algebra in case of the \schr vector field, these commutators lose their unwelcome features and enable the derivation of effective, asymptotic splittings.

In \cite{bader14eaf}, the current authors established a new framework for a numerical approach to the linear time-dependent problem with an autonomous potential,
\begin{equation}\nonumber
\partial_t u(x,t) = \ii[\varepsilon \Delta - \varepsilon^{-1} V(x) ]u(x,t),\ x\in[-1,1],\ t\geq0,
\end{equation}
where the underlying problem is considered to evolve in a certain Lie group, and the splitting of the linear operator on the right is followed by semidiscretisation. Due to the choice of a suitable Lie algebra, the authors \KKc{ were able to derive a new exponential splitting, that is the {\it asymptotic exponential splitting} of the following form:}
\begin{equation}
  \label{eq:splitting_FoCM}
  \ee^{\ii h(\ve\K-\ve^{-1} \DDD)}=
  \ee^{\frac12 W^{[0]}}\ee^{\frac12 W^{[1]}}\cdots \ee^{\frac12 W^{[s]}}\ee^{\WW^{[s+1]}} \ee^{\frac12 W^{[s]}}\cdots \ee^{\frac12 W^{[1]}}\ee^{\frac12 W^{[0]}}+\O{\ve^{2s+2}},
\end{equation}
where
\begin{Eqnarray*}
  W^{[0]}&=&W^{[0]}(h,\ve,\K,\DDD)=\O{\ve^0},\\
  W^{[k]}&=&W^{[k]}(h,\ve,\K,\DDD)=\O{\ve^{2k-2}},\qquad k=1,\ldots,s,\\
  \WW^{[k+1]}&=&\WW^{[k+1]}(h,\ve,\K,\DDD)=\O{\ve^{2s}}\PSc{.}
\end{Eqnarray*}
Here $\K$ and $\DDD$  are matrices that approximate the differential operator and multiplication by the potential $V$, respectively. \KKc{ Such asymptotic exponential splittings derived in \cite{bader14eaf}}  are superior to standard exponential splittings in a number of ways.

First of all, instead of quantifying the errors in terms of the step size, $h$, which could have been misleading due to large hidden constants, the errors are quantified in terms of the inherent semiclassical parameter $\ve$, taking into account the $\O{\ve^{-1}}$ oscillations characteristic of the semiclassical \schr equation.

Secondly, these require far fewer exponentials than  classical splittings to attain a given order. To be precise, the number of exponentials is shown to grow linearly, rather than exponentially, with the order.
Moreover, the exponents decay increasingly more rapidly in powers of $\ve$, yielding an \emph{asymptotic splitting}.

%Ultimately, even though some non-diagonal matrices will appear in the splittings, the computation of the exponentials can still be cheaply performed by low-dimensional Lanczos-methods because of their small exponents. The overall cost is cubic in the desired order, in contrast to the exponential costs of Yo\v{s}ida type splittings which becoming increasingly prohibitive once the Hamiltonian to be split features more than two terms.

%Firstly, we quantify the error not in terms of the step-size $\Delta t$ but of the small parameter $\ve$. There are three small quantities at play: $\ve, \Delta t$ and $1/M$ (where $M$ is the number of degrees of freedom in the semidiscretisation). By letting power laws govern the relationship between $\ve$ and the choices of $\Delta t$ and $M$, we express the error in the single quantity $\ve$.

%Secondly, the number of individual terms in \R{eq:splitting_FoCM} is remarkably small and it grows {\em linearly\/} with $s$ -- compare with the exponential growth, as a function of order, in the number of components of familiar splittings. The reason is that the arguments of the exponentials in \R{eq:splitting_FoCM} decay increasingly more rapidly in $\ve$.

Thirdly, each of these exponentials can be computed fairly easily. The exponents $W^{[0]}$ and $W^{[1]}$ are either diagonal or circulant matrices and their exponentials can be computed either directly or through FFT, respectively. Remaining exponents are very small and their exponentials can be computed cheaply using low-dimensional Lanczos methods.

The overall cost is quadratic in the desired order, in contrast to the exponential costs of Yo\v{s}ida type splittings which becomes increasingly prohibitive once the Hamiltonian to be split features more than two terms.

%Some of the $\RR_k$s are diagonal matrices, whereby computing the exponential is trivial. Other are circulants and can be computed with FFT. Finally, because of the small spectral radius of the arguments for sufficiently large $k$, remaining exponentials can be evaluated up to the requisite power of $\ve$ using a {\em very\/} low-dimensional Krylov subspace method. All in all, the cost of these splittings ends up being cubic in the desired order in contrast with the exponential cost of the Yo\v{s}ida method.

\KKc{The aim of the paper is to derive asymptotic exponential splittings for \schr equations with time-varying potentials. To develop such a splitting, we must first resort to the Magnus expansion.}We follow the approach of \cite{MKO} in Section~3, discretising the integrals in the Magnus expansion using Gauss-Legendre quadratures. However, unlike the traditional Magnus expansion for ODEs, we work with infinite dimensional operators to evaluate the commutators. To arrive at such a commutator-free expression, we work in the \emph{free Lie algebra} of the infinite dimensional operators $\dx^2$ and $V$ discussed in Section~2. Following the framework of \cite{bader14eaf}, a \KKc{symmetric} Zassenhaus splitting is carried out on the commutator-free Magnus expansion\KKc{, to present, eventually, the asymptotic exponential splitting of the fifth order (\ref{eq:ZM}). Obviously, following this derivation one can obtain the method of any desired order, see Table~\ref{tab.zh}. } Implementation and numerical examples are discussed in Section~4.

\KKc{Convergence and unitarity of our method follows from exactly the same argument that was presented in \cite{bader14eaf}. Namely it can be easily shown that all the exponents appearing in the derived splitting (\ref{eq:ZM}) are skew-Hermitian, hence the exponentials are unitary, which suffices for the stability of the method. Now, given consistency of our method (indeed, our scheme will be shown to be of local accuracy much higher than the order one required for the method to be consistent), we can use Lax-equivalence theorem and conclude the convergence of the method.}

\PSc{Realistic systems in quantum chemistry could involve time-dependent matrix-valued, highly oscillatory and stochastic potentials, among others. The first of these will require an extension of our Lie algebraic framework and is under active investigation, while extensions of an alternative scheme that was developed in a recent work \cite{IKS} could prove promising for oscillatory and low regularity potentials. In this approach the integrals appearing in the Magnus expansion are discretised at the very last stage, following a \KKc{symmetric} Zassenhaus splitting.}
%An alternative approach was developed in a recent work \cite{IKS}, where 

%This framework was later applied in \cite{IKS} to (\ref{eq.se}), where familiar authors derived a new asymptotic splitting. {\color {cyan} I am not sure how far we want to go here... We will be ready to conclude it later, when the paper is nearly finished. Basically we should write what was done in \cite{IKS} and why there is a need for that actually written paper}

%Once we analyse the equation (\ref{eq.se}) in a Hilbert space $\mathcal{H} = \CC{L}_2[-1,1]$,  its solution represents the probability density of finding the particle in state $x$ at time $t$, therefore  its norm equals one:
%$$
%\|u(x,t)\|^2_{L^2}=\int_{-1}^{1}|u(x,t)|^2dx=\|u(x,0)\|^2_{L^2}=1.
%$$
%We wish to preserve  this favourable feature, unitarity, both because of physical significance, and since, as will be explained in the sequel, it implies stability of the numerical method, hence convergence.
%

\section{Lie-group setting}
Following the established framework in \cite{bader14eaf}, we suppress the dependence on $x$ in (\ref{eq.se}) and analyse the following abstract ODE
\begin{equation}\label{eq:schrN}
\partial_t u(t) = \mathcal{A}(t)u(t),\ u(0)=u_0,
\end{equation}
where $\mathcal{A}(t) := \ii \varepsilon \dx^2 - \ii \varepsilon^{-1} V(t)$. Because the operator $A(t)$ belongs to $\mathfrak{u}\!\left(\mathcal{H}\right)$, the Lie algebra of (infinite-dimensional) skew-Hermitian operators acting on the Hilbert space $\mathcal{H}$, its flow is unitary and resides in $\mathcal{U}\!\left(\mathcal{H}\right)$ -- the Lie group corresponding to $\mathfrak{u}\!\left(\mathcal{H}\right)$.

The vector field in the semiclassical \schr  equation is a linear combination of the action of two operators, $\dx^2$ and multiplication by the interaction potential $V$. Since our main tools, Magnus expansion and exponential splitting methods, entail nested commutation, we consider the free Lie algebra,
\begin{displaymath}
  \GG{F}={\CC{FLA}}\{\dx^2,V\},
\end{displaymath}
i.e.,\ the  linear-space closure of all nested commutators generated by $\dx^2$ and $V$. Following \cite{bader14eaf}, we describe their action on sufficiently smooth functions, e.g.
\begin{displaymath}
  [V,\dx^2]u=V(\dx^2u)-\dx^2(Vu)=-(\dx^2V)u-2(\dx V)\dx u
\end{displaymath}
which means that $[V,\dx^2]=-(\dx^2V)-2(\dx V)\dx$. In general, we note that all terms in $\GG{F}$ belong to the set
\begin{displaymath}
  \GG{G}=\left\{\sum_{k=0}^n y_k(x)\dx^k\,:\, n\in\BB{Z}_+,\; y_0,\ldots,y_n\in {\CC{C}}_p^\infty[-1,1] \right\},
\end{displaymath}
where the subscript $p$ means periodicity in $[-1,1]$. It is trivial to observe that $\GG{G}$ is itself a Lie algebra with the commutator
\begin{align}\nonumber
	\left[ \sum_{i=0}^n  f_i(x)\partial^i_x , \sum_{j=0}^m g_j(x)\partial^j_x \right]
	=&
		\sum_{i=0}^n\sum_{j=0}^m \sum_{\ell=0}^i \binom{i}{\ell} f_i(x)  \left(\partial^{i-\ell}_x g_j(x)\right)\partial^{\ell+j}_x \\
		\label{eq:fullcommutator}
	&
		\mbox{}- \sum_{j=0}^m\sum_{i=0}^n \sum_{\ell=0}^j \binom{j}{\ell} g_j(x)  \left(\partial^{j-\ell}_x f_i(x)\right)\partial^{\ell+i}_x.
\end{align}

In similar vein to \cite{bader14eaf}, we proceed in the pursuit of stability to replace all odd powers of $\dx$ that are accompanied by $\ii$. The identities,
\begin{Eqnarray*}
  y\dx&=&-\frac12 \left[ \Int{\xi}{0}{x}{y(\xi)}\right]\dx^2 -\frac12\dx y+\frac12 \dx^2 \left[\Int{\xi}{0}{x}{y(\xi)}\,\ \cdot\,\right]\!,\\
  y\dx^3&=&-(\dx y)\dx^2-\frac14 \left[\Int{\xi}{0}{x}{y(\xi)} \right] \dx^4+\frac14 \dx^3y -\frac12 \dx^2[(\dx y)\,\cdot\,] +\frac14 \dx^4 \left[\Int{\xi}{0}{x}{y(\xi)}\,\ \cdot\,\right]\!,\\
  y\dx^5&=&\frac43 (\dx^3y)\dx^2-\frac53(\dx y)\dx^4-\frac16 \left[\Int{\xi}{0}{x}{y(\xi)}\right]\dx^6 -\frac12 \dx^5y+\frac76 \dx^2[(\dx^3y)\,\cdot\,]\\
  &&\mbox{}-\frac56 \dx^4[(\dx y)\,\cdot\,]+\frac16 \dx^6\left[\Int{\xi}{0}{x}{y(\xi)}\,\ \cdot\,\right]\!,
\end{Eqnarray*}
where $y$ is a $\CC{C}^1$ function, suffice for our presentation. The general form for expressing $y\dx^{2s+1}$ as a linear combination of even derivatives is reported in \cite{bader14eaf}.

In the Zassenhaus splitting for time-independent potentials \cite{bader14eaf}, the commutators arise solely from the {\em symmetric Baker--Campbell--Hausdorff formula} where each commutator has an odd number of {\em letters}. In the case of the \schr equation, where our operators $\dx^2$ and $V$ are each multiplied by $\ii$, this translates into an odd power of $\ii$ for each commutator.

The Magnus expansion, however, does not posses such a desirable structure -- it has commutators with odd as well as even number of letters. As a consequence, we have odd and even powers of $\ii$ accompanying our terms and it is not enough to blindly replace odd powers of $\dx$. Instead, we replace all odd powers of $\dx$ when accompanied by an odd power of $\ii$ and all even powers of $\dx$ when accompanied by an even power of $\ii$. A general formula for the replacement of even derivatives by odd derivatives can be proven along similar lines as \cite{bader14eaf}. For all practical purposes, however, we only require the identities
\begin{Eqnarray*}
  y&=&- \left[\Int{\xi}{0}{x}{y(\xi)}\right]\dx + \dx \left[\Int{\xi}{0}{x}{y(\xi)}\,\cdot\,\right]\!,\\
  y\dx^2&=&- \frac13 \left[\Int{\xi}{0}{x}{y(\xi)}\right]\dx^3  - \frac23 (\dx y) \dx - \frac13 \dx [(\dx y)\,\cdot\,] +\frac13 \dx^3 \left[\Int{\xi}{0}{x}{y(\xi)}\,\cdot\,\right]\!,\\
  y\dx^4&=&- \frac15 \left[\Int{\xi}{0}{x}{y(\xi)}\right]\dx^5  - \frac43 (\dx y)\dx^3 + \frac{8}{15}(\dx^3 y)\dx + \frac{7}{15} \dx [(\dx^3 y)\,\cdot\,] \\
  &&\mbox{}-\frac23 \dx^3 [(\dx y)\,\cdot\,] +\frac15 \dx^5 \left[\Int{\xi}{0}{x}{y(\xi)}\,\cdot\,\right]\!,
\end{Eqnarray*}
which can be easily verified directly.

%Once the appropriate odd and even derivative terms are replaced, operators of the form $f(t) \dx^k + \dx^k [f(t)\ \cdot\ ]$ start appearing ubiquitously in our workings. Far from being unique to the Magnus expansion, this is characteristic of the free Lie algebra of $\dx^2$ and $V$ and these algebraic forms also appear in Zassenhaus splittings for time-independent potentials \cite{bader14eaf}. It is worth noting that there is indeed an algebraic theory behind these structures which is under development in parallel, but not much is lost here by considering these as mere notational conveniences. We introduce a useful notation
%\[ \Ang{k}{f} := f \bullet\ \dx^k =  \frac12\left\{f \circ \dx^k + \dx^k \circ f\right\} = \frac12\left\{f \dx^k + \dx^k [f\ \cdot\ ] \right\},\qquad f \in C_p^\infty(2), \]
%where $\bullet$ is the {\em Jordan product} on the associative algebra of $\circ$ (operatorial composition). It should be obvious that $\Ang{k}{f}$ is $\O{\ve^{-k}}$, is Hermitian if $k$ is even, and skew-Hermitian if $k$ is odd. In this notation $\Ang{2}{1} = \dx^2$ and $\Ang{0}{V} = V$.

Once appropriate odd and even differential operators are replaced, operators of the form $f \dx^k + \dx^k [f\ \cdot\ ]$ start appearing ubiquitously in our analysis. Far from being unique to the Magnus expansion, they are characteristic of the free Lie algebra of $\dx^2$ and $V$ -- these algebraic forms also appear in Zassenhaus splittings for time-independent potentials \cite{bader14eaf}. We introduce a convenient notation,
\[ \Ang{k}{f} := f \bullet\ \dx^k =  \Frac12\left\{f \circ \dx^k + \dx^k \circ f\right\} = \Frac12\left\{f \dx^k + \dx^k [f\ \cdot\ ] \right\},\qquad f \in \CC{C}_p^\infty([-1,1];\BB{R}), \]
where $\bullet$ is the {\em Jordan product} on the associative algebra of $\circ$ (operatorial composition). In this notation $\Ang{2}{1} = \dx^2$ and $\Ang{0}{V} = V$.

It is worth noting that there is rich algebraic theory behind these structures which will feature in another publication, but not much is lost here by considering these as merely a notational convenience. For the purpose of this work we make observations which can be verified using the machinery of \R{eq:fullcommutator} in conjunction with the odd and even derivative replacement rules. We present identities which suffice for simplifying all commutators appearing in this work,
\begin{Eqnarray}
\label{eq:id}\left[ \Ang{4}{f}, \Ang{0}{g} \right] & =& 4 \Ang{3}{f (\dx g)} -2 \Ang{1}{3(\dx f) (\dx^2 g) + f (\dx^3 g)},\\
\nonumber\left[ \Ang{3}{f}, \Ang{0}{g} \right] & =& 3 \Ang{2}{f (\dx g)} - \Frac12 \Ang{0}{3(\dx f) (\dx^2 g) + f (\dx^3 g)},\\
\nonumber\left[ \Ang{2}{f}, \Ang{2}{g} \right] & =& 2 \Ang{3}{f (\dx g) - (\dx f) g} + \Ang{1}{2(\dx^2 f) (\dx g) - 2(\dx f) (\dx^2 g) + (\dx^3 f) g - f (\dx^3 g)},\\
\nonumber\left[ \Ang{2}{f}, \Ang{1}{g} \right] & =& \Ang{2}{2 f (\dx g) - (\dx f) g} - \Frac12 \Ang{0}{2(\dx f) (\dx^2 g) + f (\dx^3 g)},\\
\nonumber\left[ \Ang{2}{f}, \Ang{0}{g} \right] & =& 2 \Ang{1}{f (\dx g)},\\
\nonumber\left[ \Ang{1}{f}, \Ang{1}{g} \right] & =& \Ang{1}{f (\dx g)-(\dx f) g},\\
\nonumber\left[ \Ang{1}{f}, \Ang{0}{g} \right] & =& \Ang{0}{f (\dx g)}.
\end{Eqnarray}
The terms $\ii \dx^2 = \ii \Ang{2}{1}$ and $\ii V = \ii \Ang{0}{V}$ reside in
\[ \mathfrak{H} = \{ \ii^{k+1} \Ang{k}{f}\ :\ f \in \CC{C}_p^{\infty}([-1,1];\BB{R}),\ k\geq 0 \} \]
and, as evident through a few examples in \R{eq:id}, all commutators of elements of $\mathfrak{H}$ also reside in $\mathfrak{H}$. In other words, $\mathfrak{H}$ is a Lie algebra such that
\[ {\CC{FLA}}\{\ii \dx^2,\ii V\} \subseteq \mathfrak{H}, \]
and it suffices to work directly in $\mathfrak{H}$ using the rules \R{eq:id} instead of proceeding via \R{eq:fullcommutator} followed by the odd-even derivative replacement rules.

For a real valued $f$, $\Ang{k}{f}$ is symmetric if $k$ is even and skew-symmetric otherwise. This property is preserved under discretisation once we use spectral collocation on a uniform grid. In that case $\dx$ is discretised as a skew-symmetric matrix $\K$ and $V$ is discretised as a diagonal matrix $\DDD_V$. The term $\Ang{k}{f}$ is discretised as $\left(\K^k \DDD_{f} + \DDD_{f} \K^k \right)/2$ which is clearly symmetric when $k$ is even and skew-symmetric otherwise. Consequently, elements of $\mathfrak{H}$ such as $\ii^{k+1} \Ang{k}{f}$, which are skew-Hermitian operators, discretise to skew-Hermitian matrices of the form $\ii^{k+1} \left(\K^k \DDD_{f} + \DDD_{f} \K^k \right)/2$.

This structural property of $\mathfrak{H}$ is responsible for unitary evolution and numerical stability of our schemes since exponentials of skew-Hermitian matrices are unitary.

%This table already hints at the stability of the overall algorithm: Starting from terms of the form $\ii\Ang{2}{1},\ii\Ang{0}{V(t)}$, we will operate using commutators only which serves as an induction argument. As can be seen from the table, the Jordan products form a Lie algebra (closed under commutation).
%We can assume skew-Hermiticity of the left-hand sides by replacing all terms with even indices $\Ang{2k}{f}$ by $\ii\Ang{2k}{f}$\marginpar{\textcolor{red}{This sounds confusing}} and since commutators of skew-Hermitian operators are skew-Hermitian, we ensure that all generated terms will be skew-Hermitian themselves.

\newtheorem{defn}{Definition}
\begin{defn}
\label{eq:defht}
The height of a term is defined as
\[\mathrm{ht}\!\left(\sum_{i=0}^n \alpha_i \Ang{k_i}{f_i} \right) = \max\{k_1,\ldots,k_n\}. \]
\end{defn}
These terms benefit from a remarkable property of \emph{height reduction} which is stated here without proof,
\[\mathrm{ht}\left(\left[\Ang{k}{f},\Ang{l}{g} \right]\right) \leq  k+ l - 1.\]
For the commutators relevant to this work, this property can be verified by a quick inspection of the identities \R{eq:id}.

For the largest part, our work will proceed in the language of the undiscretised operators introduced in this section. At the very last stage we will resort to spectral collocation on the uniform grid over $[-1,1]$ for spatial discretisation. For this purpose we will need at least $M = \O{\ve^{-1}}$ points since (regardless of initial conditions) the solution of the \schr equation develops spatial oscillations of order $\O{\ve^{-1}}$ \cite{jin11mac,bao02ots}. Consequently, $\K$ scales like $\O{\ve^{-1}}$ and $\left(\K^k \DDD_{f} + \DDD_{f} \K^k \right)/2 = \O{\ve^{-k}}$. Keeping eventual discretisation in mind, we abuse notation and write $\Ang{k}{f} = \O{\ve^{-k}}$.

\PSc{More formally, following \cite{bao02ots} we assume that the solution $u(t)$, which is known to feature $\O{\ve^{-1}}$ oscillations, obeys the bounds,
\begin{equation}
    \left\| \partial_t^m \dx^k u(x,t) \right\| \leq C_{m,k}\, \ve^{-m-k} ,\qquad t \in [0,T].
\end{equation}
In this context
\[ \| \Ang{k}{f} u \| = \Frac12 \| \left(\K^k \DDD_{f} + \DDD_{f} \K^k \right) u \| \leq  \|V\|_{\infty} C_{k}\,\ve^{-k}  = \O{\ve^{-k}}. \]
Although it is possible to work in a more rigourous language throughout, the shorthand $\Ang{k}{f} = \O{\ve^{-k}}$ is indeed seen to be based on firm theoretical grounds while simplifying exposition greatly. We also remind the reader that the growth of derivatives of the potential, while certainly effecting error constants in our splittings (and therefore of concern in the context of moderately small values of $\ve$), are irrelevant in the asymptotic limit of $\ve \rightarrow 0$ since they don't scale with $\ve$ and don't effect the asymptotic analysis carried out here.}

The property of height reduction leads to a systematic decrease in the size of terms with commutation,
\[ \left[\Ang{k}{f},\Ang{l}{g} \right] = \O{\ve^{-k-l+1}}. \]
Going further, we want to analyse all terms in the common {\em currency} of the inherent semiclassical parameter $\ve$ and assume that our choice of the time-step, $h$, is governed by $h=\O{\ve^{\sigma}}$, for some $0<\sigma \leq 1$. Larger values of $\sigma$ correspond to very small time steps and are best avoided.

\section{The solution}

\subsection{The Magnus expansion}
To look for the solution of (\ref{eq:schrN}) one needs to take into account some features of  the operator $\mathcal{A}(t)$. First of all it depends on time and it cannot be assumed that its values in different points of time commute, i.e. we assume that $[\mathcal {A}(t_1),\mathcal{A}(t_2)]\neq0$ and give up the hope that the solution is of the simple form $\ee^{\int_0^t\A(\xi)d\xi}u_0$. Secondly $\A(t)$ evolves in a Lie algebra so the solution of (\ref{eq:schrN}) resides in a corresponding Lie group.
Both properties can be dealt with elegantly using the famous result from \cite{magnus54ote} by writing the solution as single exponential,
\begin{equation}\label{eq.MagFlow}
u(t)=\ee^{\Theta(t)}u(0),
\end{equation}
where the infinite series $\Theta(t)=\sum_{k=1}^\infty\Theta_k(t)$, also called as {\it Magnus expansion}, is an element of the underlying Lie algebra. Its convergence has been shown in \cite{iserles99ots}, \cite{moan08},  \cite{hochbruck03omi} for sufficiently small time--steps. Obviously we truncate this series and advance with adequately small time step $h$
\begin{equation}\label{eq.MagFlow_th}
u(t+h) = \ee^{\Theta(t+h,t)}u(t),
\end{equation}
starting from the initial step,
\begin{equation}\label{eq.MagFlow_th_0}
u(h) = \ee^{\Theta(h,0)}u(0),
\end{equation}
where we understand that the operator $\ee^{\Theta(t+h ,t)}$ is a flow evolving the solution from $t$ to $t+h$. \KKc{Let us observe now, that the aim of the paper consists in a derivation of the asymptotic exponential splitting for  a certain function of type $\ee^{\Theta(t+h,t)}$.  This means, that the algorithm we are going to present will advance in small time steps $h$ (exactly like the method (\ref{eq:splitting_FoCM}) does). For the clarity of \PSc{exposition}, however, we will focus on the first step of Magnus expansion}, i.e. (\ref{eq.MagFlow_th_0}), \PSc{noting that \R{eq.MagFlow_th}, when required for any time window $[t,t+h]$, is easily recovered from \R{eq.MagFlow_th_0} by a straightforward translation of the vector field $\mathcal{A}(\xi)$ to $\mathcal{A}(t+\xi)$.}
\PSc{For convenience we} shorten the notation, writing  $\Theta(h)$ instead of $\Theta(h,0)$.

Simple differentiation of the \emph{ansatz} in \eqref{eq.MagFlow_th_0} together with elementary algebra, see \cite{iserles99ots} or \cite{blanes09tme} for details, lead to the conclusion that the exponent $\Theta(t)$ satisfies the \emph{dexpinv equation},
\begin{equation}\label{1:eq:dexpinv}
	\dot{\Theta}(h)
	= \dexp^{-1}_{\Theta(h)}\A(h)
	= \sum_{k=0}^\infty \frac{B_k}{k!} \ad_{\Theta(h)}^k\A(h)
	, \qquad \Theta(0)=0,
\end{equation}
%The second identity is standard textbook knowledge \cite{blanes09tme} and
where $B_k$ are Bernoulli numbers ($B_0=1,\ B_1=-\frac{1}{2},\ B_2=\frac{1}{6}Ő\ B_3=0,\ B_4=-\frac{1}{30},\ B_5=0,\ B_6=\frac{1}{42}$) and the adjoint representation is defined recursively by $\ad^0_A V = V$ and $\ad^{k+1}_A V=[A,\ad^k_A V]$.
The solution of (\ref{1:eq:dexpinv}) is an infinite series and can be obtained using Picard iterations. It was proposed  in  \cite{magnus54ote} and widely analysed in \cite{iserles99ots,iserles00lgm,blanes09tme}.

The first few terms of the Magnus expansion ordered by size in $h$ are
\PSc{\begin{equation}\label{eq:SymMag}
\begin{split}
\Theta(h)=&\int_0^h\A(\xi)d\xi-\frac{1}{2}\int_0^h\int_0^{\xi_1}[\A(\xi_2),\A(\xi_1)]d\mathbf{\xi}\\
&+\frac{1}{12}\int_0^h\int_0^{\xi_1}\int_0^{\xi_1}[\A(\xi_2),[\A(\xi_3),\A(\xi_1)]]d\mathbf{\xi}\\
&+\frac{1}{4}\int_0^h\int_0^{\xi_1}\int_0^{\xi_2}[[\A(\xi_3),\A(\xi_2)],\A(\xi_1)]d\mathbf{\xi} + \cdots.\\
%&-\frac{1}{24}\int_0^h\int_0^{\xi_1}\int_0^{\xi_1}\int_0^{\xi_3}[\A(\xi_2),[[\A(\xi_4),\A(\xi_3)],\A(\xi_1)]]d\mathbf{\xi}\\
%&-\frac{1}{24}\int_0^h\int_0^{\xi_1}\int_0^{\xi_2}\int_0^{\xi_2}[[\A(\xi_3),[\A(\xi_4),\A(\xi_2)]],\A(\xi_1)]d\mathbf{\xi}\\
%&-\frac{1}{8}\int_0^h\int_0^{\xi_1}\int_0^{\xi_2}\int_0^{\xi_3}[[[[\A(\xi_4),\A(\xi_3)],\A(\xi_2)],\A(\xi_1)]d\mathbf{\xi}\\
%&+\cdots.
\end{split}
\end{equation}}
We say that a multivariate integral of a nested commutator, $\mathcal{I}$, is of grade $m$ if $\mathcal{I}=\mathcal{O}(h^m)$ for every smooth $\mathcal{A}$.
Truncating the Magnus expansion at grade $p$ to $\Omega_p(h)=\Theta(h)+\mathcal{O}(h^{p+1})$, preserves time symmetry \cite{iserles00lgm},  \cite{iserles01tsa}.
Time symmetry means that not only the exact flow $\varphi$, but also the numerical flow $\Phi = \ee^{\Omega_p(h)}$, satisfy
\begin{equation}
\varphi(h,0)\circ \varphi(0,h) = \Id, \qquad
\Phi_{h,0}\circ \Phi_{0,h} = \Id.
\end{equation}
As one can observe, the time symmetry of the numerical flow is equivalent to the fact that
\begin{equation}\label{eq:exp_sym}
\Omega_p(0,h)=-\Omega_p(h,0).
\end{equation}
Time symmetry is a desirable feature because truncation by power with odd $p$ leads to a gain of an extra unit of order, see \cite{iserles00lgm}. This means that if we aim for a numerical method of order six it suffices to consider the truncation of the Magnus expansion only to the terms listed in (\ref{eq:SymMag}).

\subsection{Magnus expansion in practice}
It turns out that the multivariate integrals can be efficiently computed using simple univariate quadrature rules of \citeasnoun{MKO}. We will follow their approach and evaluate the potential at the Gauss--Legendre quadrature points ($t_1=\frac{1}{2}-\frac{\sqrt{15}}{10}$, $t_2=\frac{1}{2}$, $t_3=\frac{1}{2}+\frac{\sqrt{15}}{10}$) which is then transformed \cite{iserles99ots} to obtain a far less costly quadrature. As a result, to obtain order six approximation, all the effort of approximation of the solution boils down to the following formula
\begin{equation}
\begin{split}
\label{eq:mko}
\Theta(h)=&B_1+\frac{1}{12}B_3-\frac{1}{12}[B_1,B_2]+\frac{1}{240}[B_2,B_3]+\frac{1}{360}[B_1,[B_1,B_3]]\\
&-\frac{1}{240}[B_2,[B_1,B_2]]+\frac{1}{720}[B_1,[B_1,[B_1,B_2]]]+\O{h^7},
\end{split}
\end{equation}
where
\begin{equation}\label{eq:substitutionrules}
B_1 = h\A(t_2), \ B_2=\frac{\sqrt{15}}{3} h(\A(t_3)-\A(t_1)), \ B_3=\frac{10}{3}h(\A(t_3)-2\A(t_2)+\A(t_1)).
\end{equation}
See \cite{iserles00lgm} and \cite{blanes09tme} for comprehensive information and ways to approximate the Magnus expansion using different quadrature rules and to higher orders. The former could be relevant if the time-dependent potential is only known at certain grid-points as might be the case in some control setups.
%\begin{equation}
%\dummy{V}_{0} = V(t_2), \ \dummy{V}_{1}=\frac{\sqrt{15}}{3h} (V(t_3)-V(t_1)), \ \dummy{V}_2=\frac{10}{3h^2}(V(t_3)-2V(t_2)+\A(t_1)).
%\end{equation}

Substituting $\A(t)$ with the given Hamiltonian as $\A(t) = -\ii H(t)/\ve$ and working in the free Lie algebra $\GG{H}$, we can derive a commutator free expansion using the identities \R{eq:id}.
Keeping the notation of the previous section in mind, we approximate the time derivatives of the potential by central differences, cf. \eqref{eq:substitutionrules},
\[
\dummy{V}_{0} = V(t_2), \ \dummy{V}_{1}=\frac{\sqrt{15}}{3h} (V(t_3)-V(t_1)), \ \dummy{V}_2=\frac{10}{3h^2}(V(t_3)-2V(t_2)+V(t_1)),
\]
so that
\[
B_1 = \ii h \ve \dx^2 - \ii h \ve^{-1} \dummy{V}_0, \quad
B_2 = \ii h^2 \ve^{-1} \dummy{V}_1, \quad
B_3 = \ii h^3 \ve^{-1} \dummy{V}_2.
\]
Once these are substituted in \R{eq:mko}, we use the identities \R{eq:id} along with the observation that $\dx^2 = \Ang{2}{1}$ and $V_j = \Ang{0}{V_j}$ to arrive at a Magnus expansion in the format $\sum_k \ii^{k+1} c_k  \Ang{k}{f_k}$ with $c_k \in \BB{Q}$ and $f_k \in \CC{C}_p^\infty([-1,1];\BB{R})$.

\PSc{The grade one commutators of the self-adjoint basis appearing in \R{eq:mko}, for instance, can be simplified as follows,
\begin{Eqnarray}
\nonumber [B_1,B_2] & = &  \left[ \ii h \ve \ang{1}{2} - \ii h \ve^{-1} \ang{\dummy{V}_0}{0}, -\ii h^2 \ve^{-1} \ang{\dummy{V}_1}{0} \right] \\
\label{eq:SAB_first}&=& h^3 \left[\ang{1}{2}, \ang{\dummy{V}_1}{0}\right] = 2 h^3 \ang{\dx \dummy{V}_1}{1},\\
\nonumber [B_1,B_3] & = & \left[ \ii h \ve \ang{1}{2} - \ii h \ve^{-1} \ang{\dummy{V}_0}{0}, -\ii h^3 \ve^{-1} \ang{\dummy{V}_2}{0} \right] \\
&=& h^4 \left[\ang{1}{2}, \ang{\dummy{V}_2}{0}\right] = 2 h^4 \ang{\dx \dummy{V}_2}{1},\\
\nonumber [B_2,B_3] & = & \left[ -\ii h^2 \ve^{-1} \ang{\dummy{V}_1}{0}, \ii h^3 \ve^{-1} \ang{\dummy{V}_2}{0} \right] \\
&=& 0.
\end{Eqnarray}
Consequently, the grade two commutators appearing in \R{eq:mko} are,
\begin{Eqnarray}
\nonumber [B_1,[B_1,B_3]] & = & \left[ \ii h \ve \ang{1}{2} - \ii h \ve^{-1} \ang{\dummy{V}_0}{0}, 2 h^4 \ang{\dx \dummy{V}_2}{1} \right]\\
\nonumber & = & 2 \ii h^5 \ve \left[\ang{1}{2}, \ang{\dx \dummy{V}_2}{1}\right] + 2 \ii h^5 \ve^{-1} \left[\ang{\dx \dummy{V}_2}{1}, \ang{\dummy{V}_0}{0}\right]\\
\nonumber &=& 2 \ii h^5 \ve \left(2 \ang{\dx^2 \dummy{V}_2}{2} - \Frac12 \ang{\dx^4 \dummy{V}_2}{0}\right)  + 2 \ii h^5 \ve^{-1} \ang{(\dx \dummy{V}_2)(\dx \dummy{V}_0)}{0}\\
&=& 4 \ii h^5 \ve \ang{\dx^2 \dummy{V}_2}{2} - \ii h^5 \ve \ang{\dx^4 \dummy{V}_2}{0} + 2 \ii h^5 \ve^{-1} \ang{(\dx \dummy{V}_2)(\dx \dummy{V}_0)}{0},\\
\nonumber & & \\
\nonumber [B_2,[B_1,B_2]] & = & \left[ -\ii h^2 \ve^{-1} \ang{\dummy{V}_1}{0}, 2 h^3 \ang{\dx \dummy{V}_1}{1} \right] \\
\nonumber &=& 2 \ii h^5 \ve^{-1} \left[\ang{\dx \dummy{V}_1}{1}, \ang{\dummy{V}_1}{0} \right] \\
&=& 2 \ii h^5 \ve^{-1} \ang{(\dx \dummy{V}_1)^2}{0},\\
\nonumber & & \\
\nonumber [B_1,[B_1,B_2]] &=& \left[ \ii h \ve \ang{1}{2} - \ii h \ve^{-1} \ang{\dummy{V}_0}{0}, 2 h^3 \ang{\dx \dummy{V}_1}{1} \right]\\
\nonumber & = & 2 \ii h^4 \ve \left[ \ang{1}{2},  \ang{\dx \dummy{V}_1}{1}\right] + 2 \ii h^4 \ve^{-1} \left[\ang{\dx \dummy{V}_1}{1}, \ang{\dummy{V}_0}{0}\right]\\
\nonumber &=& 2 \ii h^4 \ve \left(2 \ang{\dx^2 \dummy{V}_1}{2} - \Frac12 \ang{\dx^4 \dummy{V}_1}{0}\right) + 2 \ii h^4 \ve^{-1} \ang{(\dx \dummy{V}_1)(\dx \dummy{V}_0)}{0}\\
&=& 4 \ii h^4 \ve \ang{\dx^2 \dummy{V}_1}{2} - \ii h^4 \ve \ang{\dx^4 \dummy{V}_1}{0} + 2 \ii h^4 \ve^{-1} \ang{(\dx \dummy{V}_1)(\dx \dummy{V}_0)}{0}.
\end{Eqnarray}
The only grade three commutator that we need is
\begin{Eqnarray}
\label{eq:SAB_last}
[B_1,[B_1,[B_1,B_2]]] &=& \left[ \ii h \ve \ang{1}{2} - \ii h \ve^{-1} \ang{\dummy{V}_0}{0}, \right.\\
\nonumber & & \qquad  \qquad \left. 4 \ii h^4 \ve \ang{\dx^2 \dummy{V}_1}{2} - \ii h^4 \ve \ang{\dx^4 \dummy{V}_1}{0} + 2 \ii h^4 \ve^{-1} \ang{(\dx \dummy{V}_1)(\dx \dummy{V}_0)}{0} \right]\\
\nonumber & = & -4 h^5 \ve^2 \left[ \ang{1}{2}, \ang{\dx^2 \dummy{V}_1}{2}\right] + h^5 \ve^2 \left[\ang{1}{2}, \ang{\dx^4 \dummy{V}_1}{0}\right]\\
\nonumber & & - 2 h^5 \left[\ang{1}{2}, \ang{(\dx \dummy{V}_1)(\dx \dummy{V}_0)}{0}\right] - 4 h^5 \left[\ang{\dx^2 \dummy{V}_1}{2} ,\ang{\dummy{V}_0}{0}\right]\\
\nonumber &=& -4 h^5 \ve^2 \left(2 \ang{\dx^3 \dummy{V}_1}{3} -\ang{\dx^5 \dummy{V}_1}{1}\right) + 2 h^5 \ve^2 \ang{\dx^5 \dummy{V}_1}{1}\\
\nonumber & & - 4 h^5 \ang{(\dx^2 \dummy{V}_1)(\dx \dummy{V}_0) + (\dx \dummy{V}_1)(\dx^2 \dummy{V}_0)}{1} - 8 h^5 \ang{(\dx^2 \dummy{V}_1)(\dx \dummy{V}_0)}{1}\\
\nonumber &=& -8 h^5 \ve^2 \ang{\dx^3 \dummy{V}_1}{3} + 3 h^5 \ve^2 \ang{\dx^5 \dummy{V}_1}{1} - h^5 \ang{12(\dx^2 \dummy{V}_1)(\dx \dummy{V}_0) + 4 (\dx \dummy{V}_1)(\dx^2 \dummy{V}_0)}{1}.
\end{Eqnarray}}

\PSc{Substituting (\ref{eq:SAB_first}--\ref{eq:SAB_last}) in \R{eq:mko} gives us a truncated Magnus expansion for the \schr equation \R{eq.se} in the Lie algebra $\GG{H}$,
\begin{Eqnarray}
\nonumber \Omega_5&=& \ii h \ve \ang{1}{2} - \ii h \ve^{-1} \ang{\dummy{V}_0}{0} - \frac{1}{12} \ii h^3 \ve^{-1} \ang{\dummy{V}_2}{0}  -\frac16 h^3 \ang{\dx \dummy{V}_1}{1}\\
\nonumber & & + \frac{1}{360}\left(4 \ii h^5 \ve \ang{\dx^2 \dummy{V}_2}{2} - \ii h^5 \ve \ang{\dx^4 \dummy{V}_2}{0} + 2 \ii h^5 \ve^{-1} \ang{(\dx \dummy{V}_2)(\dx \dummy{V}_0)}{0} \right)\\
\nonumber & &-\frac{1}{120} \ii h^5 \ve^{-1} \ang{(\dx \dummy{V}_1)^2}{0} + \frac{1}{720} \left(-8 h^5 \ve^2 \ang{\dx^3 \dummy{V}_1}{3} + 3 h^5 \ve^2 \ang{\dx^5 \dummy{V}_1}{1} \right.\\
\nonumber & & \qquad \qquad \qquad \qquad \qquad \qquad \left. - h^5 \ang{12(\dx^2 \dummy{V}_1)(\dx \dummy{V}_0) + 4 (\dx \dummy{V}_1)(\dx^2 \dummy{V}_0)}{1} \right)\\
\nonumber &&\\
\label{eq:Th5MKOsigma} &=& \overbrace{\ii  h \ve \dx^2 -\ii  h \ve^{-1} \dummy{V}_0}^{\O{\ve^{\sigma-1}}} - \overbrace{\frac{1}{12}\ii  h^{3} \ve^{-1} \dummy{V}_{2} - \frac{1}{6} h^{3} \Ang{1}{ \dx \dummy{V}_1 }}^{\O{\ve^{3\sigma-1}}}\\
\nonumber && +\overbrace{\frac{1}{360}\ii  h^{5} \ve^{-1} \Big( 2 (\dx \dummy{V}_2) (\dx \dummy{V}_0) - 3 (\dx \dummy{V}_1)^{2} \Big)}^{\O{\ve^{5\sigma-1}}} -\overbrace{\frac{1}{180} h^{5} \Ang{1}{ (\dx \dummy{V}_1) (\dx^{2}\dummy{V}_0) + 3 (\dx \dummy{V}_0) (\dx^{2}\dummy{V}_1)}}^{\O{\ve^{5\sigma-1}}}\\
\nonumber &&+ \overbrace{\frac{1}{90}\ii  h^{5} \ve \Ang{2}{ \dx^{2}\dummy{V}_2 }- \frac{1}{90} h^{5} \ve^{2} \Ang{3}{ \dx^{3}\dummy{V}_1 }}^{\O{\ve^{5\sigma-1}}} - \overbrace{\frac{1}{360} \ii h^5 \ve (\dx^4 \dummy{V}_2) + \frac{1}{240} h^5 \ve^2 \ang{\dx^5 \dummy{V}_1}{1}}^{\O{\ve^{5 \sigma +1}}}\\
\nonumber &&\\
\nonumber &=& \Theta + \O{\ve^{7\sigma-1}}.
\end{Eqnarray}
For $\sigma \leq1$, the last two terms in $\Omega_5$, which are $\O{\ve^{5\sigma+1}}$, become $\O{\ve^{7\sigma-1}}$ and can be discarded. After discarding these terms, the Magnus expansion reduces to}
\begin{Eqnarray}
%\Theta & = & \overbrace{\ii  h \ve^{-1} \dummy{V}^{(0)} +  \ii  h \ve \dx^2}^{\O{\ve^0}} + \overbrace{\frac{1}{12}\ii  h^{3} \ve^{-1} \dummy{V}_{(2)} + \frac{1}{6} h^{3} \Ang{1}{ \dx^{1}\dummy{V}^{(1)} }}^{\O{\ve^2}}\\
%&& +\overbrace{\frac{1}{360}\ii  h^{5} \ve^{-1} \left( 2 (\dx^{1}\dummy{V}^{(2)}) (\dx^{1}\dummy{V}^{(0)}) - 3 (\dx^{1}\dummy{V}^{(1)})^{2} \right)}^{\O{\ve^4}}\\
%&& -\overbrace{\frac{1}{180} h^{5} \Ang{1}{ (\dx^{1}\dummy{V}^{(1)}) (\dx^{2}\dummy{V}^{(0)}) + 3 (\dx^{1}\dummy{V}^{(0)}) (\dx^{2}\dummy{V}^{(1)})}}^{\O{\ve^4}}\\
%&&- \overbrace{\frac{1}{90}\ii  h^{5} \ve \Ang{2}{ \dx^{2}\dummy{V}^{(2)} }+ \frac{1}{90} h^{5} \ve^{2} \Ang{3}{ \dx^{3}\dummy{V}^{(1)} }}^{\O{\ve^4}}  + \O{\ve^6}.
\label{eq:Th5MKO1} \Omega_5 & = & \overbrace{\ii  h \ve \dx^2 -\ii  h \ve^{-1} \dummy{V}_0}^{\O{\ve^{\sigma-1}}} - \overbrace{\frac{1}{12}\ii  h^{3} \ve^{-1} \dummy{V}_{2} - \frac{1}{6} h^{3} \Ang{1}{ \dx \dummy{V}_1 }}^{\O{\ve^{3\sigma-1}}}\\
\nonumber && +\overbrace{\frac{1}{360}\ii  h^{5} \ve^{-1} \Big( 2 (\dx \dummy{V}_2) (\dx \dummy{V}_0) - 3 (\dx \dummy{V}_1)^{2} \Big)}^{\O{\ve^{5\sigma-1}}}\\
\nonumber && -\overbrace{\frac{1}{180} h^{5} \Ang{1}{ (\dx \dummy{V}_1) (\dx^{2}\dummy{V}_0) + 3 (\dx \dummy{V}_0) (\dx^{2}\dummy{V}_1)}}^{\O{\ve^{5\sigma-1}}}\\
\nonumber &&+ \overbrace{\frac{1}{90}\ii  h^{5} \ve \Ang{2}{ \dx^{2}\dummy{V}_2 }- \frac{1}{90} h^{5} \ve^{2} \Ang{3}{ \dx^{3}\dummy{V}_1 }}^{\O{\ve^{5\sigma-1}}} = \Theta + \O{\ve^{7\sigma-1}}.
\end{Eqnarray}
%\begin{Eqnarray*}
%%\Theta & = & \overbrace{\ii  h \ve^{-1} \dummy{V}^{(0)} +  \ii  h \ve \dx^2}^{\O{\ve^0}} + \overbrace{\frac{1}{12}\ii  h^{3} \ve^{-1} \dummy{V}_{(2)} + \frac{1}{6} h^{3} \Ang{1}{ \dx^{1}\dummy{V}^{(1)} }}^{\O{\ve^2}}\\
%%&& +\overbrace{\frac{1}{360}\ii  h^{5} \ve^{-1} \left( 2 (\dx^{1}\dummy{V}^{(2)}) (\dx^{1}\dummy{V}^{(0)}) - 3 (\dx^{1}\dummy{V}^{(1)})^{2} \right)}^{\O{\ve^4}}\\
%%&& -\overbrace{\frac{1}{180} h^{5} \Ang{1}{ (\dx^{1}\dummy{V}^{(1)}) (\dx^{2}\dummy{V}^{(0)}) + 3 (\dx^{1}\dummy{V}^{(0)}) (\dx^{2}\dummy{V}^{(1)})}}^{\O{\ve^4}}\\
%%&&- \overbrace{\frac{1}{90}\ii  h^{5} \ve \Ang{2}{ \dx^{2}\dummy{V}^{(2)} }+ \frac{1}{90} h^{5} \ve^{2} \Ang{3}{ \dx^{3}\dummy{V}^{(1)} }}^{\O{\ve^4}}  + \O{\ve^6}.
%\Omega_5 & = & \overbrace{\ii  h \ve \dx^2 -\ii  h \ve^{-1} \dummy{V}_0}^{\O{\ve^0}} - \overbrace{\frac{1}{12}\ii  h^{3} \ve^{-1} \dummy{V}_{2} - \frac{1}{6} h^{3} \Ang{1}{ \dx \dummy{V}_1 }}^{\O{\ve^2}}\\
%&& +\overbrace{\frac{1}{360}\ii  h^{5} \ve^{-1} \Big( 2 (\dx \dummy{V}_2) (\dx \dummy{V}_0) - 3 (\dx \dummy{V}_1)^{2} \Big)}^{\O{\ve^4}}\\
%&& -\overbrace{\frac{1}{180} h^{5} \Ang{1}{ (\dx \dummy{V}_1) (\dx^{2}\dummy{V}_0) + 3 (\dx \dummy{V}_0) (\dx^{2}\dummy{V}_1)}}^{\O{\ve^4}}\\
%&&+ \overbrace{\frac{1}{90}\ii  h^{5} \ve \Ang{2}{ \dx^{2}\dummy{V}_2 }- \frac{1}{90} h^{5} \ve^{2} \Ang{3}{ \dx^{3}\dummy{V}_1 }}^{\O{\ve^4}} = \Theta + \O{\ve^6}.
%\end{Eqnarray*}
\PSc{We note that, due to the property of height reduction discussed in Section~2, a grade $n$ commutator in the Magnus expansion of $\mathcal{A}(t)$ should be $\O{\ve^{n \sigma -1}}$. This can indeed be verified in the above expansion. Asymptotically speaking, in terms of $\ve$, the terms in the expansion are decreasing in size with increasing $n$ for any $\sigma>0$, so that convergence of the Magnus expansion also occurs for much larger time steps such as $h=\O{\ve^{1/2}}$ or $h=\O{\ve^{1/4}}$. This is a considerable improvement over existing analysis.}

Since $\Omega_5$ includes the term $\ii  h \ve \dx^2-\ii  h \ve^{-1} \dummy{V}_0$, its exponential is, at the very least, as troublesome to approximate as the problem of solving the \schr equation with time-independent potential. Fortunately the Zassenhaus procedure is sufficiently flexible and can tackle such modified Hamiltonians with ease.

\subsection{Zassenhaus}
Let us recall the basic principle for the iterative symmetric Zassenhaus splitting \cite{bader14eaf}.
Our goal is to compute $e^{\mathcal{W}^{[0]}}$, where $\mathcal{W}^{[0]}=X+Y$ and $X,Y = \O{\ve^{p}}$.
Using the \emph{symmetric Baker-Campbell-Hausdorff} (sBCH) formula \cite{dynkin47eot,casas09aea}, we then write
 %$$e^{-\frac12 hA}e^{h(A+B)}e^{-\frac12 hA}=e^{\sBCH(-hA,h(A+B))},$$
\begin{equation}\label{eq:2:sbch_inverted}
  e^{\mathcal{W}^{[0]}} =e^{\frac12 X} e^{\sBCH(-X,\mathcal{W}^{[0]})} e^{\frac12 X}.
\end{equation}
Grade-three commutators of $X$ and $Y$ are at most $\O{\ve^{3p}}$ so that $\sBCH(-X,X+Y) = Y + \O{\ve^{3p}}$. Thus we have extracted $X$ from the exponent at the cost of correction terms in form of higher-order commutators. Assuming that the corrections are decreasing in size, it is then enough to identify the largest terms as $W^{[1]}$ in the central exponent $\mathcal{W}^{[1]}=\sBCH(-X,\mathcal{W}^{[0]})$ and to continue the iteration until the desired accuracy is reached,
\def\WW{\mathcal{W}}
\begin{equation} \label{eq:definitionWW}
	\WW^{[k+1]} = \sBCH(-W^{[k]}, \WW^{[k]}), \quad \WW^{[0]}=X+Y.
\end{equation}
In this notation, the splitting after $s$ steps can be written as
\begin{align*}
	\exp(X+Y) &= e^{\frac12 W^{[0]}}e^{\frac12 W^{[1]}}
									\cdots e^{\frac12 W^{[s]}} e^{\WW^{[s+1]}} e^{\frac12 W^{[s]}}
									\cdots e^{\frac12 W^{[1]}}e^{\frac12 W^{[0]}}.
\end{align*}
We emphasize that, in principle, we can freely choose the elements $W^{[k]}$ that we want to extract. Except for some special cases, at least one of the exponents in this splitting will feature an infinite series of terms. To construct a finite splitting scheme featuring a certain accuracy we may discard, at each stage, all terms smaller than the desired threshold.

Assuming that a grade $k$ commutator of $X$ and $Y$ scales as $\O{\ve^{k p }}$, convergence of the series requires $p>0$ at the very least. In the case of the \schr equation, we choose $X = \ii  h \ve \dx^2$, $p=\sigma-1$, and this  naively translates to a very stringent time step restriction: $\sigma > 1$. However, the remarkable feature of height reduction means that a grade $k$ commutator in this context scales as $\O{\ve^{k p + (k - 1)}}$ and convergence requirements become significantly milder: we need $p + 1> 0$ which translates to $\sigma > 0$.

\subsection{Zassenhaus on Magnus}
We perform a Zassenhaus splitting on $\Omega_5$, choosing to extract the largest terms -- analysed in powers of $\ve$ -- first. We commence the splitting with $W^{[0]} = \ii  h \ve \dx^2$, although we could equally well choose  $W^{[0]} = -\ii  h \ve^{-1} \dummy{V}_0$, for instance, and arrive at a variant of the splitting presented here. The exponent to be split is $\WW^{[0]} = \Omega_5$ and the first step involves computing the sBCH formula. Here, once again, the rules of the free Lie algebra $\GG{H}$, \R{eq:id}, suffice for arriving at a commutator free expression,
\begin{align*}
	\WW^{[1]} &= \sBCH(-W^{[0]}, \WW^{[0]})\\
    &= -\overbrace{\ii  h \ve^{-1} \dummy{V}_0}^{\O{\ve^{\sigma-1}}} + \overbrace{\frac{1}{12}\ii  h^3 \ve^{-1} \Big(2(\dx \dummy{V}_0)^2 - \dummy{V}_2\Big) -\frac{1}{6} h^{3} \Ang{1}{ \dx \dummy{V}_1 } +\frac{1}{6}\ii  h^{3} \ve \Ang{2}{ \dx^{2} \dummy{V}_0 } }^{\O{\ve^{3\sigma-1}}} \\
& -\overbrace{\frac{1}{24}\ii  h^{3} \ve (\dx^{4} \dummy{V}_0)}^{\O{\ve^{3\sigma+1}}} - \overbrace{\frac{1}{360}\ii  h^{5} \ve^{-1} \Big(8(\dx \dummy{V}_0)^{2} (\dx^{2} \dummy{V}_0) +3(\dx \dummy{V}_1)^{2} - 12(\dx \dummy{V}_2) (\dx \dummy{V}_0)\Big)}^{\O{\ve^{5\sigma-1}}}\\
& +\overbrace{\frac{1}{30} h^{5} \Ang{1}{ 2(\dx \dummy{V}_0) (\dx^{2} \dummy{V}_1) - (\dx \dummy{V}_1) (\dx^{2} \dummy{V}_0) }}^{\O{\ve^{5\sigma-1}}}\\
& -\overbrace{\frac{1}{720}\ii  h^{5} \ve \Ang{2}{ 127(\dx \dummy{V}_0) (\dx^{3} \dummy{V}_0)+ 130(\dx^{2} \dummy{V}_0)^2 - 18(\dx^{2} \dummy{V}_2)}}^{\O{\ve^{5\sigma-1}}}\\
& +\overbrace{\frac{1}{60} h^{5} \ve^{2} \Ang{3}{ \dx^{3} \dummy{V}_1 } - \frac{13}{90}\ii  h^{5} \ve^{3} \Ang{4}{ \dx^{4}\dummy{V}_0 }  }^{\O{\ve^{5\sigma-1}}}  + \O{\ve^{7\sigma-1}}.
\end{align*}
At the second stage we select the largest remaining element $W^{[1]} = -\ii  h \ve^{-1} \dummy{V}_0$, whereby
\begin{align*}
	\WW^{[2]} &= \sBCH(-W^{[1]}, \WW^{[1]})\\
& = \overbrace{\frac{1}{12}\ii  h^{3} \ve^{-1}  \Big( 2(\dx \dummy{V}_0)^{2} - \dummy{V}_2 \Big) - \frac{1}{6} h^{3} \Ang{1}{ \dx \dummy{V}_1 } +\frac{1}{6}\ii  h^{3} \ve \Ang{2}{ \dx^{2}\dummy{V}_0 }}^{\O{\ve^{3\sigma-1}}}\\
&-\overbrace{\frac{1}{24}\ii  h^{3} \ve  (\dx^{4}\dummy{V}_0)}^{\O{\ve^{3\sigma+1}}} - \overbrace{\frac{1}{360}\ii  h^{5} \ve^{-1} \Big(13(\dx \dummy{V}_0)^{2} (\dx^{2}V) + 3(\dx \dummy{V}_1)^{2} - 12 (\dx \dummy{V}_2) (\dx \dummy{V}_0)\Big)}^{\O{\ve^{5\sigma-1}}}\\
&+\overbrace{\frac{1}{30} h^{5} \Ang{1}{ 2 (\dx \dummy{V}_0) (\dx^{2} \dummy{V}_1) - (\dx \dummy{V}_1) (\dx^{2} \dummy{V}_0) }}^{\O{\ve^{5\sigma-1}}}\\
&-\overbrace{\frac{1}{720}\ii  h^{5} \ve \Ang{2}{ 127(\dx \dummy{V}_0) (\dx^{3} \dummy{V}_0)+ 130(\dx^{2} \dummy{V}_0)^2 - 18(\dx^{2} \dummy{V}_2)}}^{\O{\ve^{5\sigma-1}}}\\
& +\overbrace{\frac{1}{60} h^{5} \ve^{2} \Ang{3}{ \dx^{3} \dummy{V}_1 } - \frac{13}{90}\ii  h^{5} \ve^{3} \Ang{4}{ \dx^{4}\dummy{V}_0 }  }^{\O{\ve^{5\sigma-1}}}  + \O{\ve^{7\sigma-1}}.
\end{align*}
We terminate the procedure by letting $W^{[2]}$ consist of the $\O{\ve^{3\sigma-1}}$ terms in $\WW^{[2]}$ and are left with $\O{\ve^{5\sigma-1}}$ and $\O{\ve^{3\sigma+1}}$ terms in $\WW^{[3]} = \WW^{[2]} - W^{[2]}$ once we ignore $\O{\ve^{7\sigma-1}}$ terms. Since $\O{\ve^{3\sigma+1}}$ terms can be subsumed into the $\O{\ve^{5\sigma-1}}$ terms for $\sigma \leq 1$, combining them in this way is not a cause for concern.  The outcome is the splitting,
\begin{equation}
\label{eq:ZM}
	\ee^{\Omega_5} = \ee^{\frac12 W^{[0]}}\ee^{\frac12 W^{[1]}} \ee^{\frac12 W^{[2]}} \ee^{\WW^{[3]}} \ee^{\frac12 W^{[2]}} \ee^{\frac12 W^{[1]}}\ee^{\frac12 W^{[0]}} + \O{\ve^{7\sigma-1}},
\end{equation}
with
\begin{eqnarray*}
%W^{[0]} & = & \ii   \ve h\dx^2 = \O{\ve^0},\\
%W^{[1]} & = & \ii   \ve^{-1}h \dummy{V}^{(0)} = \O{\ve^0},\\
%W^{[2]} & = & \frac{1}{12}\ii   \ve^{-1} h^3\left(2 (\dx \dummy{V}^{(0)})^2 + \dummy{V}^{(2)} \right) + \frac{1}{6} h^{3} \Ang{1}{ \dx \dummy{V}^{(1)} } -\frac{1}{6}\ii  \ve h^{3} \Ang{2}{ \dx^{2} \dummy{V}^{(0)} } = \O{\ve^{2}},\\
%\WW^{[3]} & = & \frac{1}{24}\ii \ve h^3 (\dx^{4} \dummy{V}^{(0)}) \\
%& & + \frac{1}{360} \ii \ve^{-1} h^5 \left( 13 (\dx \dummy{V}^{(0)})^{2} (\dx^{2} \dummy{V}^{(0)}) - 3 (\dx \dummy{V}^{(1)})^{2} + 12 (\dx \dummy{V}^{(2)}) (\dx \dummy{V}^{(0)}) \right)\\
%& & +\frac{1}{30} h^5 \Ang{1}{ 2 (\dx V) (\dx^{2}\dummy{V}^{(1)}) - (\dx \dummy{V}^{(1)}) (\dx^{2}V) } \\
%& & - \frac{1}{720} \ii \ve h^5 \Ang{2}{ 127 (\dx \dummy{V}^{(0)}) (\dx^{3} \dummy{V}^{(0)}) + 130 (\dx^{2} \dummy{V}^{(0)})^{2} + 18 (\dx^{2}\dummy{V}^{(2)})} \\
%& & -\frac{1}{60}  \ve^{2} h^5 \Ang{3}{ \dx^{3}\dummy{V}^{(1)} } + \frac{13}{90}\ii   \ve^{3} h^5 \Ang{4}{ \dx^{4} \dummy{V}^{(0)} } = \O{\ve^4}.
W^{[0]} & = & \ii   \ve h\dx^2 = \O{\ve^{\sigma-1}},\\
W^{[1]} & = & -\ii   \ve^{-1}h \dummy{V}_0 = \O{\ve^{\sigma-1}},\\
W^{[2]} & = & \frac{1}{12}\ii   \ve^{-1} h^3\Big(2 (\dx \dummy{V}_0)^2 - \dummy{V}_2 \Big) - \frac{1}{6} h^{3} \Ang{1}{ \dx \dummy{V}_1 } +\frac{1}{6}\ii  \ve h^{3} \Ang{2}{ \dx^{2} \dummy{V}_0 } = \O{\ve^{3\sigma-1}},\\
\WW^{[3]} & = & -\frac{1}{24}\ii \ve h^3 (\dx^{4} \dummy{V}_0)  - \frac{1}{360} \ii \ve^{-1} h^5 \Big( 13 (\dx \dummy{V}_0)^{2} (\dx^{2} \dummy{V}_0) + 3 (\dx \dummy{V}_1)^{2} - 12 (\dx \dummy{V}_2) (\dx \dummy{V}_0) \Big)\\
& & +\frac{1}{30} h^5 \Ang{1}{ 2 (\dx V) (\dx^{2}\dummy{V}_1) - (\dx \dummy{V}_1) (\dx^{2}V) } \\
& & - \frac{1}{720} \ii \ve h^5 \Ang{2}{ 127 (\dx \dummy{V}_0) (\dx^{3} \dummy{V}_0) + 130 (\dx^{2} \dummy{V}_0)^{2} - 18 (\dx^{2}\dummy{V}_2)} \\
& & +\frac{1}{60}  \ve^{2} h^5 \Ang{3}{ \dx^{3}\dummy{V}_1 } - \frac{13}{90}\ii   \ve^{3} h^5 \Ang{4}{ \dx^{4} \dummy{V}_0 } = \O{\ve^{5\sigma-1}}.
\end{eqnarray*}

%@article{iserles00lgm,
%	Author = {A. Iserles and H.Z. Munthe-Kaas and S.P. N{\o}rsett and A. Zanna},
%	Journal = {Acta Numerica},
%	Pages = {215--365},
%	Title = {Lie-group methods},
%	Volume = {9},
%	Year = {2000}}

%\PBc{\subsection{Derivation of the scheme at a glance}

\section{A numerical scheme}
\PBc{As we have seen, the derivation of the method has two components.
First, we choose the desired order of accuracy in the small parameter $\varepsilon$ and compute the Magnus expansion up to this order $\Omega_p$.
This will lead to an effective exponent of the form \eqref{eq:mko}, \PSc{detailed steps for which} can be found in \PSc{\cite{MKO,blanes09tme}}.
\PSc{Commencing from these expansions we compute the commutator-free Magnus expansion using the rules \R{eq:id} of the Lie algebra $\GG{H}$.}
Once we have computed this effective Hamiltonian, we start the Zassenhaus algorithm, detailed in Table~\ref{tab.zh}.
\begin{table}[tbh]
\caption{\label{tab.zh}Zassenhaus algorithm}
\begin{tabular*}{\textwidth}{l}
	\toprule\\[-3mm]	
		$\begin{array}{l}
			\textbf{Initialise:} \\[1mm]
			%W^{[0]}:=-\ii h \varepsilon^{-1} V_0,\;
			\WW^{[0]}:=\PSc{\Omega_{2p+1},\; \text{exploiting time-symmetry for a local error of $\O{h^{2p+3}}$}}\\[1mm]
			W^{[0]}:= \ii   \ve h\dx^2,\;k=0
			\\
			%W^{[1]}:= -\ii   \ve^{-1}h \dummy{V}_0 = \O{\ve^{\sigma-1}},\\
			\midrule
			\textbf{while } \PSc{k\leq p}\\
			\;\;\begin{array}{ll}
			\WW^{[k+1]}&:=\sBCH(-W^{[k]},\WW^{[k]})\\
			\WW^{[k+1]}-W^{[k+1]}&:=\O{\varepsilon^{(2k+1)\sigma-1}} \text{ implicitly defines } W^{[k+1]}\\
			k &:= k+1\\
			\end{array}\\
			\textbf{end while}\\
	\text{Final method:}\\
	\PSc{\ee^{\Omega_{2p+1}}}=\ee^{\frac12 W^{[0]}}\ee^{\frac12 W^{[1]}}\cdots \PSc{\ee^{\frac12 W^{[p]}}\ee^{W^{[p+1]}} \ee^{\frac12 W^{[p]}}}\cdots \ee^{\frac12 W^{[1]}}\ee^{\frac12 W^{[0]}}+\O{\ve^{(2p+3)\sigma-1}}
	\end{array}$\\
  \bottomrule
\end{tabular*}
\end{table}}

For numerical realisation of these splittings schemes it is typical to impose periodic boundary conditions in order to resolve spatial oscillations with spectral accuracy. Recall that we restrict the domain to $[-1,1]$, imposing periodic boundaries at $x=\pm 1$. \PSc{We discretise using spectral collocation on the equispaced grid $\{x_n\}_{n=-N}^N$, $x_n = n/(N+\frac12)$}, $|n|\leq N$, where $M=2N+1$ is the number of grid points. The unknowns are $u_n\approx u(x_n)$, $|n|\leq N$. The differential operator $\dx$ is discretised as a circulant matrix $\K$ and $V$ as a diagonal $\DDD_V$.

All exponents in our splitting \R{eq:ZM} are of the form $\ii^{k+1} \Ang{k}{f}$ and are discretised as skew-Hermitian matrices,
\[ \PSc{\ii^{k+1} \Ang{k}{f} \leadsto \ii^{k+1} (\DDD_f \K^k + \K^k \DDD_f)/2.} \]
Since the exponential of a skew-Hermitian matrix is unitary, unitary evolution and (consequently) unconditional stability of the method are guaranteed.

The outermost exponentials $W^{[0]}$ and $W^{[1]}$ are replaced by the circulant $\ii \ve h\K^2$ and the diagonal matrix $-\ii \ve^{-1}h \DDD_{\dummy{V}_0}$, respectively.
\PSc{The lowest order scheme of the type \R{eq:ZM} can be obtained by ignoring the exponents $W^{[2]}$ and $\WW^{[3]}$ from \R{eq:ZM},
\[ \MM{u}^{n+1} = \ee^{\frac12 \ii h \ve \K^2} \ee^{-\ii h \ve \DDD_{\dummy{V}_0}} \ee^{\frac12 \ii h \ve \K^2} \MM{u}^{n}, \]
which features an $\O{\ve^{3\sigma-1}}$ error since the largest term ignored is $W^{[2]}$.
Clearly, this is the Strang splitting after freezing the potential in the middle of the interval. }

The exponential of the circulant matrix \PSc{$\frac12 \ii \ve h\K^2$} is evaluated to machine precision using Fast Fourier Transform (FFT) in $\O{M \log M}$ operations while the diagonal matrix is exponentiated directly in $\O{M}$ operations.

\PSc{The first non-trivial splitting is obtained upon including $W^{[2]}$ once more,
\[ \MM{u}^{n+1} = \ee^{\frac12 \ii h \ve \K^2} \ee^{-\frac12 \ii h \ve \DDD_{\dummy{V}_0}} \ee^{\tilde{W}^{[2]}} \ee^{-\frac12 \ii h \ve \DDD_{\dummy{V}_0}} \ee^{\frac12 \ii h \ve \K^2} \MM{u}^{n}, \]
where
\begin{Eqnarray*}
\tilde{W}^{[2]} & =&\frac{1}{12}\ii   \ve^{-1} h^3\DDD_{2 (\dx \dummy{V}_0)^2 - \dummy{V}_2 } - \frac{1}{12} h^{3} \left(\DDD_{ \dx \dummy{V}_1 } \K + \K \DDD_{ \dx \dummy{V}_1 } \right) \\
&&\quad+\frac{1}{12}\ii  \ve h^{3} \left(\DDD_{ \dx \dummy{V}_1 } \K^2 + \K^2 \DDD_{\dx^{2} \dummy{V}_0 } \right),
\end{Eqnarray*}
is the discretised version of $W^{[2]}$. This splitting commits an error of $\O{\ve^{5\sigma-1}}$. We remind the reader that if the pursuit of a splitting with an $\O{\ve^{5\sigma-1}}$ error was the objective, it would suffice to start with $\Omega_3$, which is a lower order Magnus truncation that is easier to obtain and uses merely two Gauss--Legendre quadrature knots, while to obtain splittings that are higher order than the $\O{\ve^{7\sigma-1}}$ splitting given in \R{eq:ZM} we would need to commence with a higher order Magnus expansion than those discussed here.}

\PSc{The exponents $W^{[2]}$ and $\WW^{[3]}$ appearing in the non-trivial splittings} do not posses a structure amenable to exact exponentiation. However, they are very small -- $\O{\ve^{3 \sigma -1}}$ and $\O{\ve^{5 \sigma -1}}$, respectively. For $\sigma = 1$, the most costly case we consider, the exponentials of these terms can be evaluated to $\O{\ve^6}$ accuracy using merely three and two Lanczos iterations, respectively \PSc{\cite{bader14eaf}.
These iterations involve the computation of $\tilde{W}^{[2]} \MM{v}$, which can be achieved using a few FFTs (remember that $\K^k$ in $\DDD_{ \dx \dummy{V}_1 } \K + \K \DDD_{ \dx \dummy{V}_1} $ and $\DDD_{ \dx \dummy{V}_1 } \K^2 + \K^2 \DDD_{\dx^{2} \dummy{V}_0}$ is a circulant).}

We refer the curious reader to \cite{bader14eaf,IKS} where semi-discretisation strategies, stability analysis and exponentiation methods are addressed in greater detail.

\subsection{A numerical example}
\begin{figure}[h]
\begin{center}
    \includegraphics[trim = 0.4cm 0.5cm 0.4cm 0cm, clip = true, width=184pt]{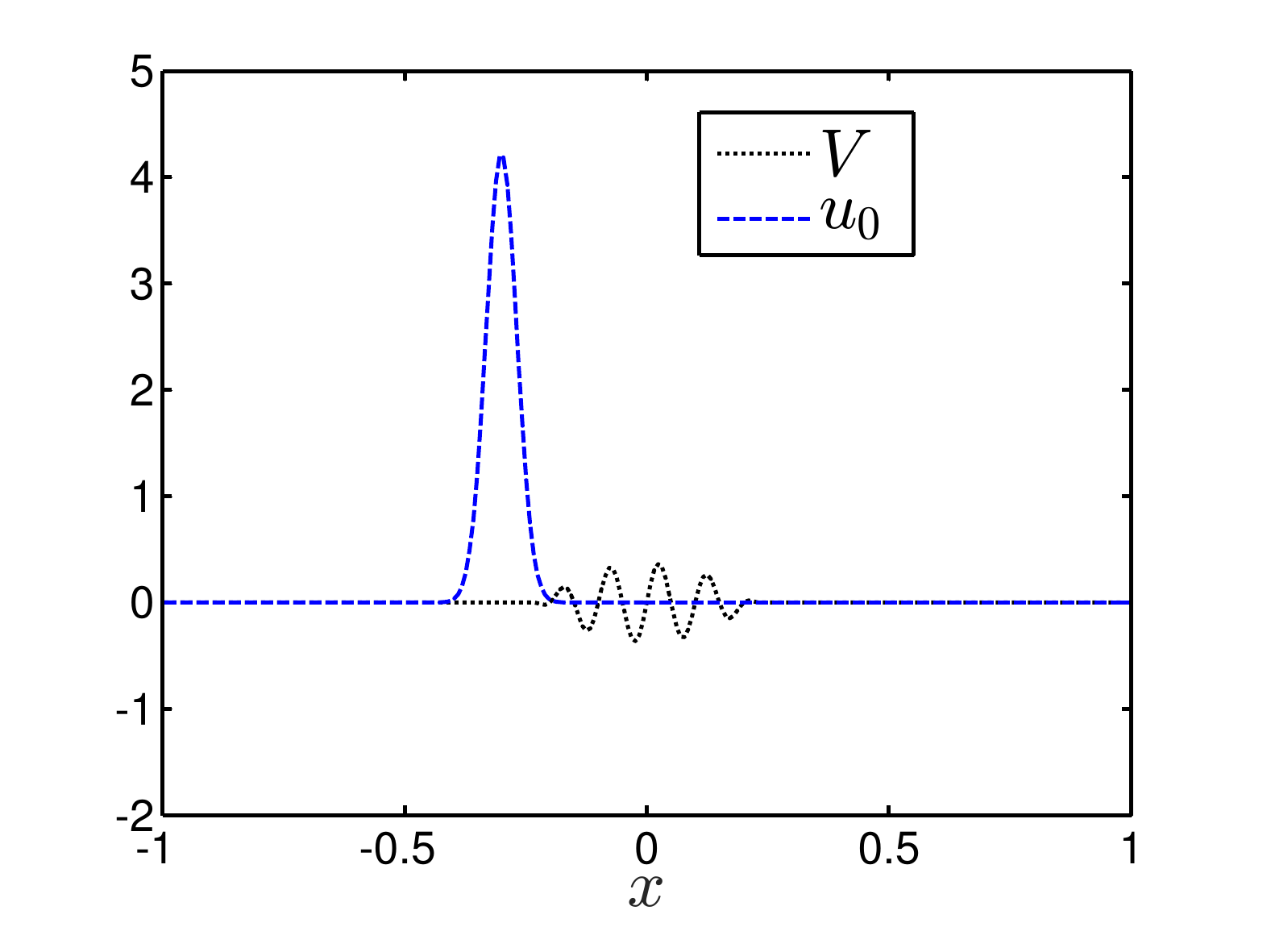}
    \includegraphics[trim = 0.2cm 0.5cm 0.6cm 0cm, clip = true, width=184pt]{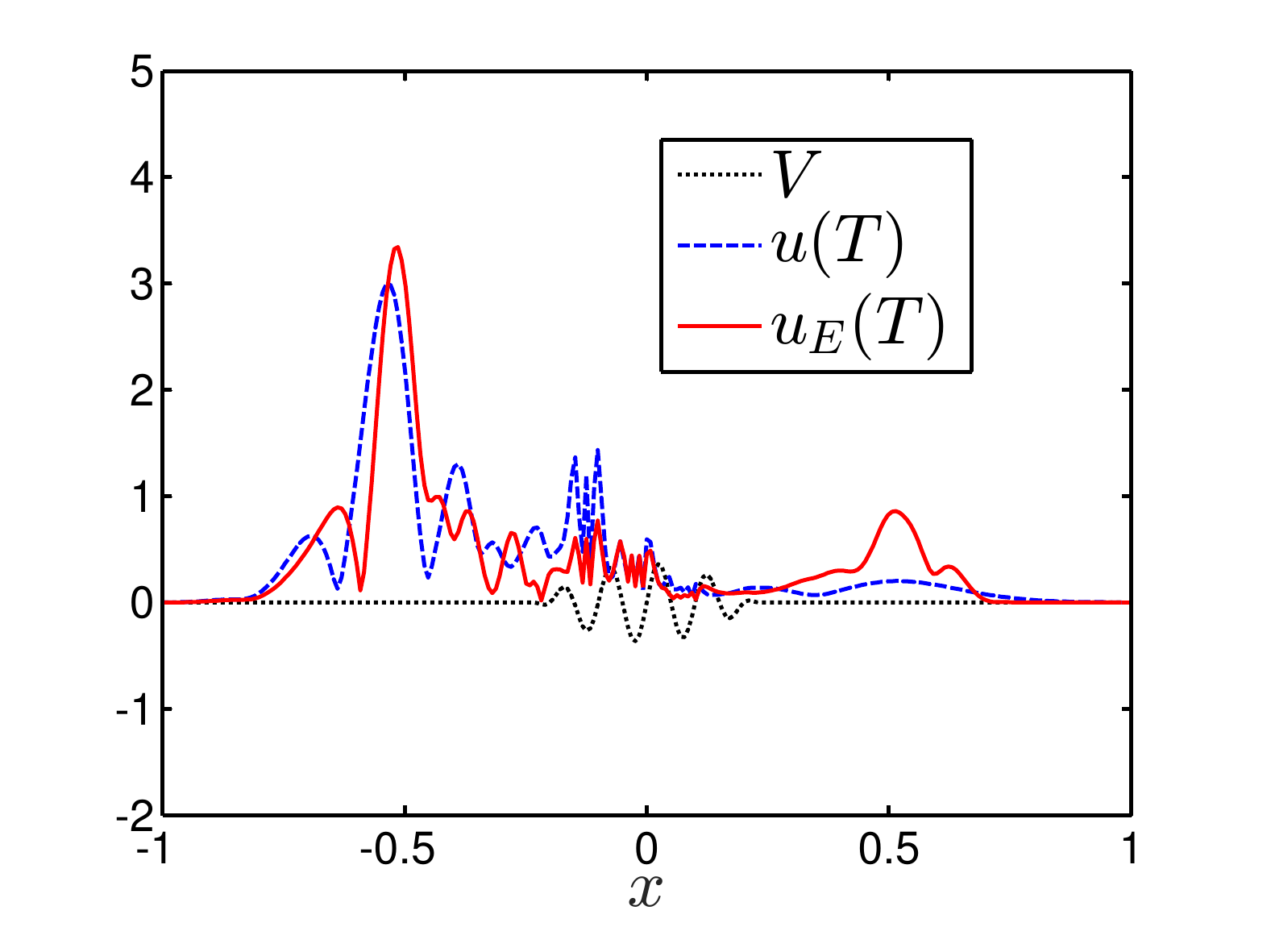}
    \caption{ (left) initial wave-packet $u_0$; (right) final wave-packets at time $T=0.75$: $u(T)$ under the influence of $V_0$ and $u_E(T)$ under the influence of $V_E(x,t) = V_0 + E(x,t)$. }
    \label{fig:DW}
\end{center}
\end{figure}

Consider the evolution of the wave-packet
\begin{equation*}
\label{eq:u0}
u_0(x) = (\delta \pi)^{-1/4} \exp\left(\ii k_0 \frac{(x-x_0)}{\delta}-\frac{(x-x_0)^2}{2 \delta }\right)
\end{equation*}
with $x_0 = -0.3$, $k_0=0.1$ and $\delta = 1.22\times 10^{-4}$, heading towards the lattice potential
\begin{equation*}
    \label{eq:V0}
    V_0 = \rho(4 x)\sin(20\pi x),
\end{equation*}
where
\[ \rho(x) = \begin{cases} \exp\left(-1/(1-x^2)\right) &\mbox{for } |x| < 1,\\
0 &\mbox{otherwise, } \end{cases}\]
is a bump function. When the semiclassical parameter is $\ve = 2^{-8}$, the wave-packet evolves to $u(T)$ at $T=0.75$ (Figure~\ref{fig:DW}) under the influence of the time-independent potential $V_0$ alone. When we excite it using an additional time-varying potential,
\begin{equation*}
\label{eq:L}
E(x,t) = \rho(3t-1)\rho(\sin(2\pi(x-t))),
\end{equation*}
so that the wave packet evolves under $V_E(x,t) = V_0(x) + E(x,t)$, a significantly larger part of the wave packet is able to make it across the lattice to the right hand side (see $u_E(T)$ in Figure~\ref{fig:DW}).

\begin{figure}[h]
\begin{center}
    \includegraphics[trim = 0.4cm 0.5cm 0.4cm 0cm, clip = true, width=184pt]{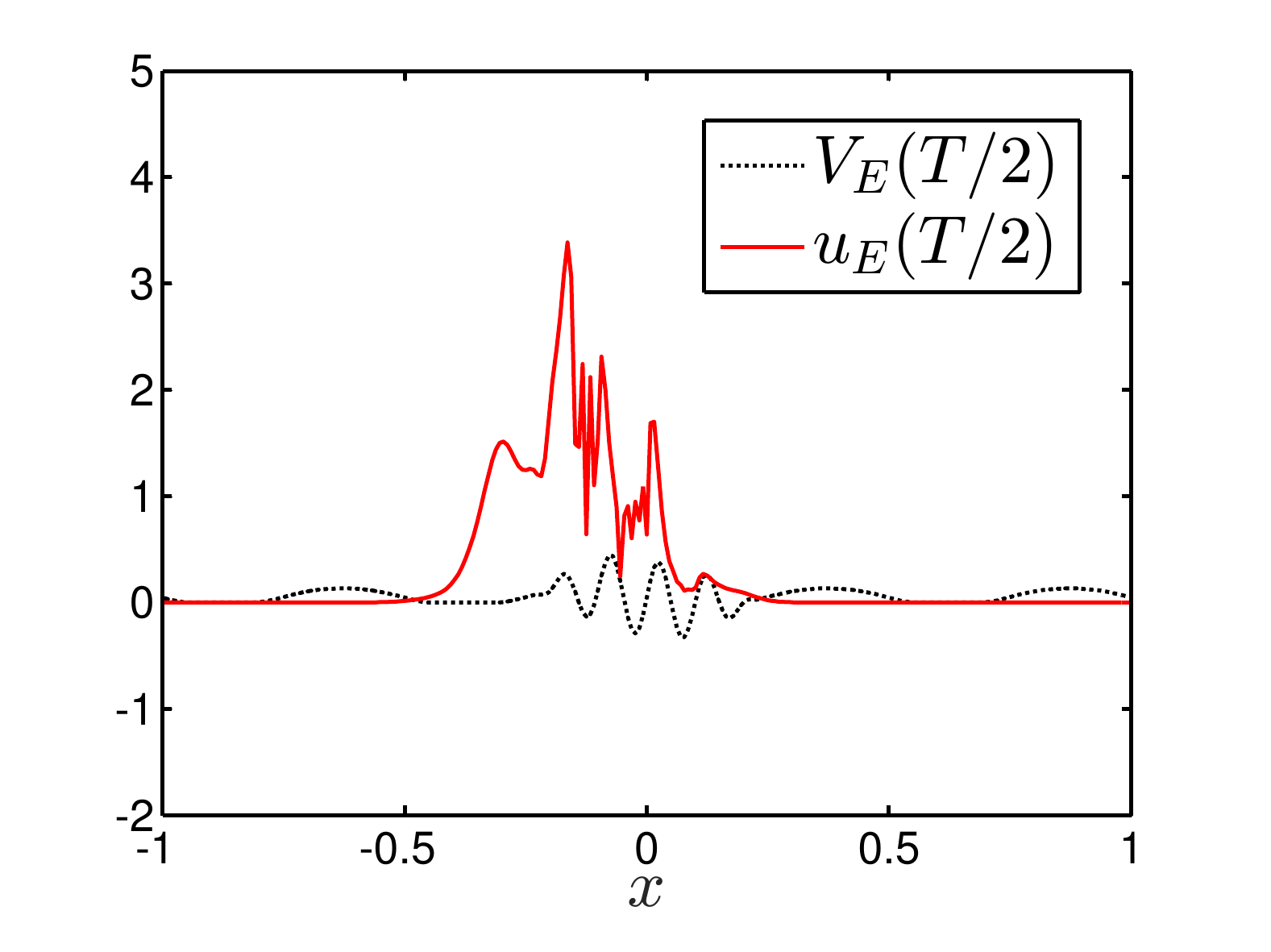}
    \includegraphics[trim = 0.2cm 0.5cm 0.6cm 0cm, clip = true, width=184pt]{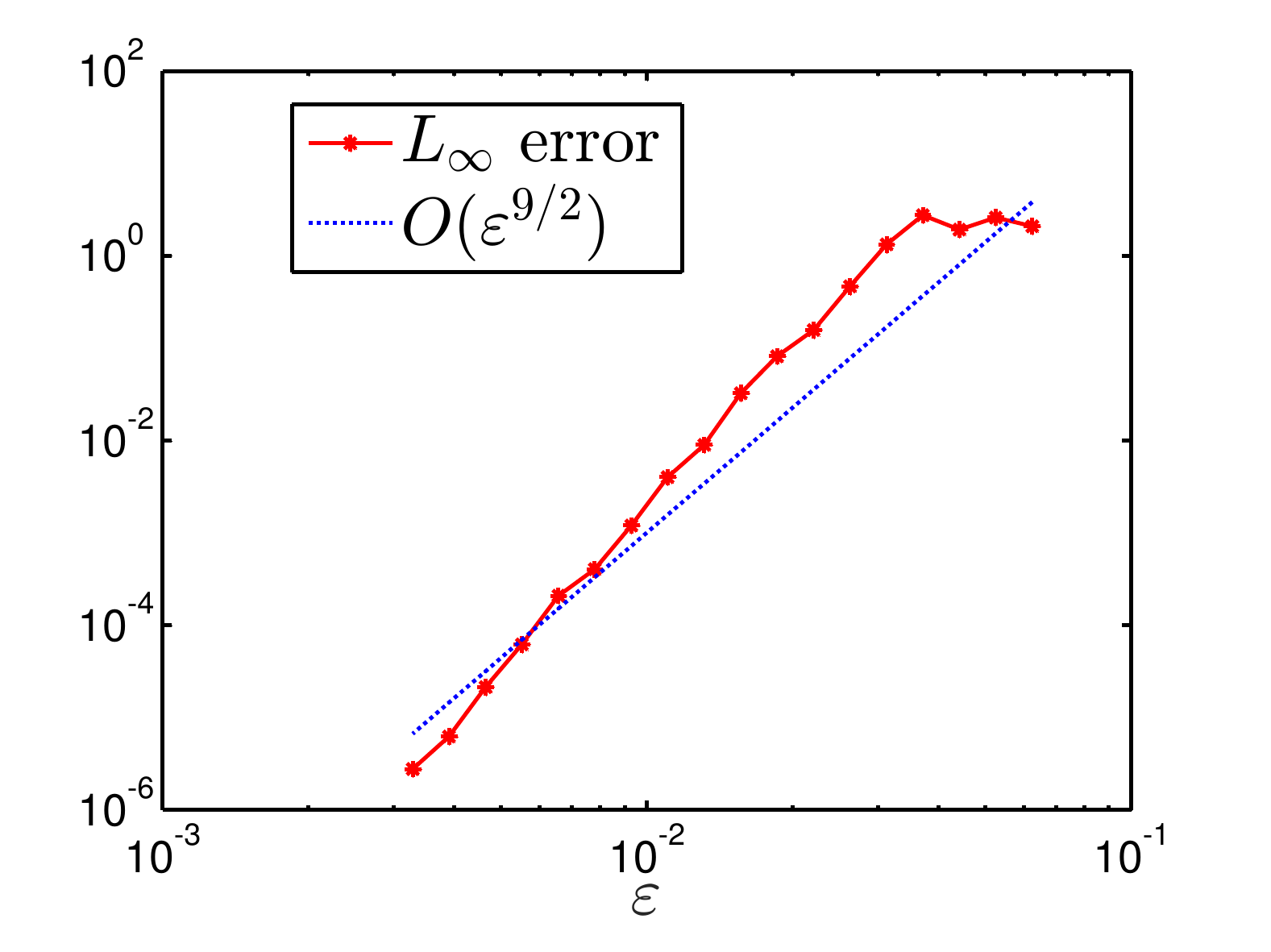}
    \caption{(left) effective potential at the middle of the time interval, $t=T/2$ (where $T=0.75$);
    (right) global error of the Zassenhaus splitting \R{eq:ZM} at $T=0.5$.}
    \label{fig:Tby2}
\end{center}
\end{figure}

The excitation pulse is not active for the entire duration since $\rho(3t-1)$ acts as a smooth envelope simulating the switching on and off of the time-varying component of the potential. The excited potential is evident at $t=T/2$ in Figure~\ref{fig:Tby2}.

In Figure~\ref{fig:Tby2} we present the global error at time $T=0.5$ in the propagation of $u_0$ to $u_E$ under the influence of $V_E$ using the scheme \R{eq:ZM}.
Under the scaling $\sigma =1$, we commit a local $\CC{L}_2$ error of $\O{\ve^4}$ per time step in the splitting scheme \R{eq:ZM}. Since the number of time steps is $\O{\ve^{-\sigma}}$, the global error is $\O{\ve^5}$. The precise scaling used in our experiments is $M \sim 5 \ve^{-1}$ and $h \sim 2 \ve$.

Our analysis has been in the context of the $\CC{L}_2$ inner product and the corresponding norm which, upon discretisation, translates to an $\ell_2$ norm scaled by a factor of $\sqrt{2/M}$. Where $\CC{L}_\infty$ error is of greater interest, it should be noted that $\| \MM{v} \|_{\ell_\infty} \leq \sqrt{M/2} \| \MM{v} \|_{\ell_2}$ and consequently we may expect the global $\CC{L}_\infty$ error to be $\O{\ve^{9/2}}$ for $\sigma = 1$. This is indeed seen to be the case through numerical experiments in Figure~\ref{fig:Tby2}.

\subsection{Finding a reference solution}
\PSc{Since no analytic solution of \R{eq.se} is available, reference solutions must also be obtained via a numerical approach. We obtain the reference solution $\MM{u}_R$ for our numerical experiments by resorting to a Strang splitting,
\[ \MM{u}_R(t+h_R) = \ee^{\Frac12 h_R \ve \K^2}\ \ee^{-h_R \ve^{-1} \DDD_{V(t+h_R/2)}}\ \ee^{\Frac12 h_R \ve \K^2}\ \MM{u}_R(t), \]
where $\MM{u}_R \in \BB{C}^{M_R}$ lives on a much finer grid than the solution of \R{eq:ZM}. In each of the $T/h_R$ time steps required for finding the solution $\MM{u}_R(T)$, the potential is frozen in the middle of the interval $[t,t+h_R]$.}

\PSc{Since such a splitting is also the lowest order in the Magnus--Zassenhaus family of schemes, we require very small times steps for convergence -- certainly $h_R \ll h$ is required for the reference solution to possess an error smaller than the scheme \R{eq:ZM} whose error we are attempting to quantify.}

\PSc{We rely on this method for producing reliable reference solutions since it is simple and its error is easily analysed. Directly exponentiating a Hamiltonian (via MATLAB's {\tt expm}, for instance) with potential frozen at the middle of the interval is more expensive but no more accurate than the Strang splitting -- this is because freezing the potential is akin to disregarding the nested integrals and commutators in the Magnus expansion which are of the same size (in powers of $\ve$) as the error committed in the Strang splitting.}

\PSc{Another factor we must take into account is the growth of spatial oscillations with decreasing $\ve$. To capture this, starting from $M_R = 3 M = 15 \ve^{-1}$, we iteratively increase the grid resolution for the reference solution till no high frequencies are clipped and convergence is achieved. In the end, the spatial resolution used for obtaining a reference solution is much greater than that used for \R{eq:ZM}, $M_R \gg M \sim 5 \ve^{-1}$.}

%take larger and larger constants of proportionality while choosing $M_R=\O{\ve^{-1}}$ (the degrees of freedom for spatial discretisation) for computing the reference solution. The error in the reference solution is estimated by generating another reference solution with a larger number of grid points, $3 M_R$, and smaller time steps, $h_R/2$, to which the candidate reference solution is compared. A reference solution is only accepted when this difference is very small and much smaller than the error of the scheme \R{eq:ZM}.

\PSc{Using such a low order method for generating reference solutions to a high degree of accuracy in a brute force manner means generating reference solutions is orders of magnitude slower than the splitting method \R{eq:ZM} requiring validation. The exorbitant cost of reference solutions is what restricts experimental study of numerical errors to moderate values of $\ve$ and $T$.}

\bibliographystyle{agsm}
%{agsm}\begin{eqnarray}
\bibliography{Schrodinger_BIKS}

\end{document}